\documentclass{article}

\usepackage{arxiv}

\usepackage[utf8]{inputenc} 
\usepackage[T1]{fontenc}    
\usepackage{hyperref}       
\usepackage{url}            
\usepackage{booktabs}       
\usepackage{amsfonts}       
\usepackage{nicefrac}       
\usepackage{microtype}      
\usepackage{lipsum}
\usepackage{graphicx}
\graphicspath{ {./images/} }

\usepackage[figuresright]{rotating}

\usepackage{url}
\usepackage{algorithmic}
\usepackage{algorithm}
\usepackage{xcolor}

\usepackage{amsmath,amsthm}
\newtheorem{theorem}{Theorem}
\newtheorem{proposition}{Proposition}
\newtheorem{lemma}{Lemma}

\theoremstyle{definition}

\newtheorem{remark}{Remark}
\usepackage{amssymb}
\newcommand{\R}{\mathbb{R}}

\hypersetup{colorlinks=true,
	linkcolor=blue,
	urlcolor=blue,
	citecolor=red} 
\theoremstyle{remark}
\newtheorem{definition}{Definition}

\title{A preconditioner for solving linear programming problems with dense columns}

\author{
Catalina J. Villalba
\\
  Instituto de Matemática, Estatística e Computação Científica\\
Universidade Estadual de Campinas (UNICAMP)\\
Campinas, Brazil\\
  \texttt{catajv@ime.unicamp.br} \\
   \And
 Aurelio R.L. Oliveira \\
 Instituto de Matemática, Estatística e Computação Científica\\
Universidade Estadual de Campinas (UNICAMP)\\
Campinas, Brazil\\
  \texttt{aurelio@ime.unicamp.br} 
}

\begin{document}
\maketitle
\begin{abstract}
The Interior-Point Methods are a class for solving linear programming problems that rely upon the solution of linear systems. At each iteration, it becomes important to determine how to solve these linear systems when the constraint matrix of the linear programming problem includes dense columns. In this paper, we propose a preconditioner to handle linear programming problems with dense columns, and we prove theoretically that the final linear system to solve is uniformly bounded when the Interior-Point Method is converging to an optimal solution. This result is illustrated through computational experiments, which show that our proposed method is robust and competitive in terms of running time and/or number of iterations compared with existing methods.    
\end{abstract}


\section{Introduction}

The Interior-Point Methods have become an alternative approach for solving linear programming problems. In comparison to the Simplex method, one of its key advantages is the ability to obtain an optimal solution after a very low number of iterations. However, each iteration can be computationally expensive, given that it requires solving at least one linear system  (\cite{gondzio25years, vanderbei2020linear, wright1997primal}). 

The predictor-corrector version presented by Mehrotra in (\cite{mehrotra1992implementation}) is currently one of the most used methods. This approach requires the solution of two linear systems, each having the same matrix, however with a different right-hand side vector. There are two main possible approaches to solve the linear systems, the first one is the normal equations, where the Cholesky factorization can be used to solve it. The other possibility is the augmented form of the system, with an indefinite matrix, which is sparser than the normal equations, however, it has a larger dimension. 

In certain situations, the constraint matrix may include at least one dense column, indicating a substantial presence of non-zero entries. When considering both dense and sparse columns collectively, there is a risk of generating normal equations with nearly full matrices. This could result in a significant number of floating-point operations or even potentially lead to memory challenges when dealing with large-scale problems.

Some strategies have been developed by different authors to address these situations. Some have proposed modifications to the constraint matrix, leading to the derivation of equivalent problems (\cite{gondzio, vanderbei}). Alternatively, other authors have suggested methods such as adjusting the Cholesky factorization through the implementation of backward and forward substitutions (\cite{Choleskyproduct}) or solving linear systems without resorting to the typical normal equations (\cite{andersen, choi1990further}). Direct methods can also be employed using the Sherman-Morrison-Woodbury formula (\cite{wright1997primal}).

This paper introduces a preconditioner designed for a modified augmented system to address linear programming problems with dense columns. The application of the preconditioner results in a system that transforms the linear equations into a positive definite matrix. We demonstrate that this matrix remains uniformly bounded during the convergence of Interior-Point Methods toward a solution.

From a computational standpoint, leveraging the identified matrix properties allows for the utilization of iterative methods in solving linear systems. In our numerical experiments, we observed a remarkably low number of iterations required to attain a solution.

The article is structured as follows: Section \ref{LPsection} introduces some important results regarding the interior-point method. Section \ref{sec3} introduces the new preconditioner, including the numerical stability analysis and theoretical properties of the resulting linear system. Section \ref{implem}. Section \ref{experim} includes the computational experiments, and covers important aspects considered for the implementation and the density criterion. The paper concludes with final remarks in Section \ref{conclusion}.

\section{Linear programming}
\label{LPsection}

The objective of a Linear Programming (LP) problem is to optimize a linear function subject to linear constraints, including both equalities and inequalities. By adding extra variables, referred to as slack and surplus variables, the \textit{standard form} transforms inequalities into equalities. It is written in matrix form
\begin{align*}
	\min &\hspace{2mm} c^Tx\\
	\text{s.t.  }& Ax=b\\
	&x\geq 0,
\end{align*}

\noindent where $A\in\R^{m\times n}$ and we assume it has full row-rank, this is, $\textit{rank}(A)=m\leq n$; $c,x\in \R^n$, and $b\in\R^m$. And $m$ represents the number of restrictions and $n$ the number of variables. The previous linear programming problem problem is known as the primal problem. Each primal problem has associated a dual problem given by
\begin{align*}
	\max & \hspace{2mm}b^Ty\\
	\text{s.t.  }& A^Ty+z= c\\
	&y \text{ free}, z\geq 0,
\end{align*}
where $z\in\R^m$ is a slack dual variable. 

Some linear programming problems include free and upper bounded variables as follows 
\begin{align*}
	\min &\hspace{2mm} c^Tx\\
	\text{s.t.  }& Ax=b\\
	&x_i\geq 0, \hspace{2mm}i\in \mathcal{N},\\
    &0\leq x_i \leq u_i, \hspace{2mm} i\in \mathcal{U},\\
    &x_i \mbox{ free}, \hspace{2mm} i\in \mathcal{F},
\end{align*}
with $\mathcal{N} \cup \mathcal{U} \cup \mathcal{F}=\{ 1,\dots, n \}$  and they are disjoint index sets. To get a linear programming problem without free variables, we define for $i\in \mathcal{F}$:
\begin{equation*}
    x_i = x_i^+ - x_i^-,
\end{equation*}
where $x_i^+, x_i^-\geq 0$. Therefore, the above problem becomes a linear programming problem with non-negative and upper-bounded variables
\begin{equation}
\label{primal}
    \begin{split}
        \min &\hspace{2mm} c^Tx\\
    	\text{s.t.  }& Ax=b\\
        & x_i +v_i= u_i, \hspace{2mm} i\in \mathcal{U}.\\
        &x \geq 0, v_i\geq 0, \hspace{2mm}i\in \mathcal{U},\\
    \end{split}
\end{equation}
	
\noindent with $v$ a slack variable. In order to eliminate the index associated with set $\mathcal{U}$, if $|\mathcal{U}|=l$, we define the matrix $E:l\times n$ such that it contains the $i-$th columns of identity matrix if $i\in \mathcal{U}$, this is, the primal problem \eqref{primal} is rewritten as
\begin{equation}
\label{LPgeneral}
    \begin{split}
        \min &\hspace{2mm} c^Tx\\
    	\text{s.t.  }& Ax=b\\
        & Ex +v= u,\\
        &x \geq 0, v\geq 0.\\
    \end{split}
\end{equation}

The dual form of \eqref{LPgeneral} is
\begin{equation*}
    \begin{split}
	\max &\hspace{2mm} b^Ty-u^Tw\\
	\text{s.t.  }& A^Ty+z-E^Tw=c\\
	&z\geq 0, w\geq 0, 
\end{split}
\end{equation*}
where $z\in\R^m$ is the slack variable.

For simplicity, we suppose $\mathcal{U}=\{ 1,\dots ,n\}$. Therefore, the primal problem is
\begin{equation}
\begin{split} 
\min \hspace{2mm}&c^Tx\\
\textit{s.t. }& Ax=b\\
&x+v= u\\
&x,v\geq 0,
\end{split}
\label{LPgeneralfinal}
\end{equation}
and the dual problem is given by
\begin{equation}
    \begin{split} 
    \max \hspace{2mm}&b^Ty-u^Tw\\
    \textit{s.t. }& A^Ty+z-w=c\\
    &z,w\geq 0.
    \end{split}
    \label{Dualgeneralfinal}
\end{equation}

The optimality conditions for primal-dual problems are given in the following theorem from (\cite{gondzio25years}).

\begin{theorem}[Optimality conditions]
	\label{optimal}
	$(x,v)$ and $(y,z,w)$ are optimal solutions to the primal \eqref{LPgeneralfinal} and dual \eqref{Dualgeneralfinal} problems, if and only if
	\begin{align*}
		\left. \begin{array}{r}
			Ax=b\\
         x +v= u \\
        x, v\geq 0
		\end{array}\right\rbrace& \hspace{1cm} \text{Primal feasibility}\\
		\left. \begin{array}{r}
			A^Ty+z-w=c\\
	       z,w\geq 0
		\end{array}\right\rbrace  &\hspace{1cm}\text{Dual feasibility}\\
        \left. \begin{array}{r}
			XZe=0\\
            VWe=0,
		\end{array}\right\rbrace& \hspace{1cm} \text{Complementary}
	\end{align*}

 \end{theorem}

where $X$, $Z$, $V$ and $W$ are diagonal matrices with entries of vectors $x$, $z$, $v$ and $w$ respectively; and $e$ is the vector of ones. Given that computing the exact solution might be costly, we can go on to perturb the complementary condition, such that $x_jz_j=\mu$, $j=1,\dots ,n$, and $v_iw_i=\mu$, $i=1,\dots ,m,$ where $\mu$ is a constant close to zero.

Newton's method can be employed to handle the nonlinearity of the perturbed complementarity conditions and yield linear equations. To compute the search directions $(\Delta x, \Delta v, \Delta y, \Delta z, \Delta w)$, it is necessary to solve the linear system at each $k$ iteration of Interior-Point Methods:
\begin{equation*}
	\begin{split}
		\begin{bmatrix}
			A&0&0&0&0\\
		0&A^T&0&I&-I\\
            I&0&I&0&0\\
		Z^{k}&0&0&X^{k}&0\\
            0&0&W^{k}&0&V^{k}
		\end{bmatrix}
		\begin{bmatrix}
			\Delta x^k\\
			\Delta y^k\\
			\Delta v^k\\
                \Delta z^k\\
                \Delta w^k
		\end{bmatrix}&=
  \begin{bmatrix}
		b-Ax^k\\
		c-A^Ty^k-z^k+w^k\\
            u-x^k-v^k\\
		\mu e-X^kZ^ke \\
            \mu e- V^kW^ke 
	\end{bmatrix}
 :=
		\begin{bmatrix}
			r_b^{k}\\
			r_c^{k}\\
			r_u^{k}\\
            r_{xz}^{k} \\
            r_{wv}^{k}
		\end{bmatrix}.
	\end{split}
\end{equation*}

By eliminating variables we obtain the \textit{augmented system}:
\begin{align*}
\begin{bmatrix}
    (\Lambda^k)^{-1}&A^T\\
    A&0
\end{bmatrix}
\begin{bmatrix}
    \Delta x^k\\
    \Delta y^k
\end{bmatrix}
=
\begin{bmatrix}
    r_b^k\\
    r_c^k - (X^k)^{-1}r_{xz}^k-(V^k)^{-1}W^kr_u^k+(V^k)^{-1}(W^k)^{-1}r_{wv}^k
\end{bmatrix}.
\end{align*}

where $\Lambda^k = (Z^k)^{-1} X^k+(W^k)^{-1}V^k$. And eliminating $\Delta x^k$ we obtain the \textit{normal equations system}
\begin{equation}
    \label{eqnormal}
	A\Lambda^k A^T\Delta y^k= r_b^k -A\Lambda^k \left[  (X^k)^{-1}r_{xz}^k+(V^k)^{-1}W^kr_u^k+(W^k)^{-1}r_{wv}^k-r_c^k) \right]:= q^k.
\end{equation}

One of the most used Interior-Point Methods is the Predictor-Corrector developed by Mehrotra  (\cite{mehrotra1992implementation}). To improve the convergence, we split Newton's direction $\Delta^k$ into two parts, $\Delta^k=\Delta_p^k+\Delta_c^k$, where $\Delta_p^k$ is known as the predictor step or affine direction, and $\Delta_c^k$ as the corrector step. The predictor step $\Delta_p^k  = (\Delta_p x^k , \Delta_p v^k , \Delta_p y^k , \Delta_p z^k , \Delta_p w ^k)^T$ is the solution to the linear system
\begin{equation}
	\label{affine}
		\begin{cases}
			A \Delta_p x^k=r_b^k,\\
			A^{T}\Delta_p y^k+\Delta_p z^k-\Delta_p w^k=r_c^k,\\
            \Delta_p x^k+ \Delta_p v^k = r_u^k,\\
			Z^k\Delta_p x^k+X^k\Delta_p z^k=r_{xz}^k,\\
            W^k\Delta_p v^k+V^k\Delta_p w^k=r_{wv}^k.
		\end{cases}
\end{equation}

Once we solve the affine scale system and assuming $r_b^k=0$, $r_c^k=0$ and $r_u^k=0$, we define the corrector direction $\Delta_c^k=(\Delta_c x^k, \Delta_c v^k, \Delta_c y^k, \Delta_c z^k, \Delta_c w^k)^T$ as follows
\begin{equation*}
		\begin{cases}
			A \Delta_c x^k=0,\\
			A^{T}\Delta_c y^k+\Delta_c z^k-\Delta_c w^k=0,\\
            \Delta_c x^k+ \Delta_c v^k = 0,\\
			Z^k\Delta_c x^k+X^k \Delta_c z^k=\sigma^k \mu^k e - D_x^k D_z^ke=r_{c1}^k,\\
            W^k\Delta_c v^k+V^k\Delta_s w^k=\sigma^k \mu^k e - D_w^k D_v^ke=r_{c2}^k,
		\end{cases}
\end{equation*}

\noindent where $D_x^k =\text{diag}(\Delta_p x^k )$, $D_z^k =\text{diag}(\Delta_p z^k )$, $D_v^k =\text{diag}(\Delta_p v^k )$, $D_w^k =\text{diag}(\Delta_p w^k )$,  and
\begin{align*}
\mu^k &= \frac{\langle x^{k},z^{k} \rangle + \langle w^{k},v^{k} \rangle}{2n},\\
	\sigma^k  &= \left(   \frac{\langle x^{k}+\alpha_p\Delta_c x^k ,z^{k}+\alpha_d\Delta_c z^k \rangle + 
    \langle v^{k}+\alpha_p\Delta_p v^k ,w^{k}+\alpha_d\Delta_c w^k \rangle }{\langle x^{k},z^{k} \rangle + \langle w^{k},v^{k} \rangle} \right)^{3},
\end{align*}
where $\langle \cdot , \cdot \rangle$ denotes the usual inner product.

It can be observed that both systems share the same matrix, but each one has a different right-hand side. It follows that variable elimination can be used to simplify each system, leading to the solution of the normal equations or augmented system for the predictor and corrector directions. 

Concerning the normal equations, given that $A$ is assumed to be a full-row rank matrix and we are working with interior points, the normal equation is a symmetric positive definite matrix. This property allows the application of direct methods such as Cholesky Factorization (\cite{ARVK89}) to solve the linear systems. Alternatively, one can apply preconditioners to enhance the numerical stability of the normal matrix. A notable method for this purpose is the Modified Controlled Cholesky Factorization, introduced in (\cite{bocanegra2007using}), or the preconditioners defined in (\cite{CLOC17, OS03}). Moreover, considering the ill-conditioning in the final iterations of Interior-Point Methods (\cite{stability}), hybrid approaches are viable. Strategies like those proposed by (\cite{BCO13}) advocate switching the preconditioner during the Interior-Point method process to mitigate stability issues.

Regarding the augmented matrix, it is an indefinite matrix of size $m+n$. This matrix offers the advantage of avoiding the product $A\Lambda A^T$, which may result in greater fill-in in the final matrix. Given its sparse nature, we may apply techniques described in (\cite{IterativeMethods, aasen1971reduction, bunch1971direct}) or through iterative methods like MINRES (\cite{MINRES}).

To simplify the notation, the index $k$ associated with the iteration of the Interior-point method will be eliminated.

\subsection{Dense columns}
A column is considered dense if it has more than a certain number of non-zero entries. Let $P$ be a column permutation matrix such that $AP= [S,D]$, where $D: m\times k$ contains the $k$ dense columns and $S:m\times (n-k)$ the sparse columns. Also, for the diagonal matrix $\Lambda$, we have $P\Lambda P^T = \textit{diag} (\Lambda_D, \Lambda_S)$, where $\Lambda = \textit{diag}(\Lambda_{S_1},\dots , \Lambda_{S_{n-k}}):(n-k)\times (n-k)$ is a diagonal matrix whose entries are the corresponding values $x_j/z_j+ v_j/w_j$, where $j$ indexes a column of the sparse matrix $S$, and $\Lambda=\textit{diag}(\Lambda_{D_1},\dots , \Lambda_{D_k}): k\times k$ is related to the dense matrix $D$. For simplicity, assume $P=I_n$. Therefore, the normal equation is $A\Lambda A^T = S\Lambda_S S ^T+ D\Lambda_D D^T$. 

Some studies have focused on analyzing linear programming problems with matrices exhibiting these features. Proposed solutions include the modification of the $A$ matrix by splitting each dense column into several sparse columns (\cite{gondzio}), the creation of an equivalent augmented system (\cite{vanderbei}), the adaptation of the Schur’s Complement Method (\cite{andersen}), and a modification of the Cholesky factorization of the sparse matrix $S\Lambda_S S^T$ (\cite{Choleskyproduct}).

\section{A new preconditioner}
\label{sec3}

Let $A$ be the constraint matrix of a linear programming problem, such that $A=[S,D]$, where $D:m\times k$ contains the $k$ dense columns and $S:m\times (n-k)$ the sparse columns. Define $F:m\times d$, with $d\leq k$ as a matrix where each column has one element different from zero, and such that the Cholesky factorization of $S\Lambda_S S^{T}+ FF^{T}$ is well-defined. 

Define $K:(m+k+d)\times (m+k+d)$ as follows:
\begin{equation*}
    K=
    \begin{bmatrix}
    S\Lambda_S S^{T}+FF^{T}     &  D\Lambda_D^{1/2} & F\\
    \Lambda_D^{1/2} D^{T}     & -I_k&\\
    F^T  & &I_d
    \end{bmatrix}.
\end{equation*}
The matrix $K$ is indefinite and the linear system
\begin{equation}
\label{augmentedsystem}
    K\omega =
   \begin{bmatrix}
    S\Lambda_S S^{T}+FF^{T}     &  D\Lambda_D^{1/2} & F\\
    \Lambda_D^{1/2} D^{T}     & -I_k&\\
    F^T  & &I_d
    \end{bmatrix}
    \begin{bmatrix}
    \Delta y\\ q_1\\ q_2
    \end{bmatrix} = 
    \begin{bmatrix}
    q\\  0\\0
    \end{bmatrix} :=\tilde{q}.
\end{equation}
is equivalent to \eqref{eqnormal}, this is, $\Delta y$ is the Newton direction of Interior-Point Method. $q_1:k\times 1$ and $q_2: d\times 1$. 

Let $M$ be a nonsingular matrix. Instead of focusing on system \eqref{augmentedsystem}, we may solve the equivalent system $M^{-1}KM^{-T}\hat{w} = \hat{q}$, where $\hat{\omega}=M^{T}\omega$ and $\hat{q}= M^{-1}\tilde{q}$. The preconditioned system is $M^{-1}KM^{-T}\tilde{\omega} = \tilde{q}$. Define $M^{-1}$ as follows
\begin{equation*}
M^{-1}=
    \begin{bmatrix}
        L^{-1}&L^{-1}D\Lambda_D^{1/2}&-L^{-1}F\\
        0&\Lambda_D^{-1/2}&0\\
        0&0&I_k
    \end{bmatrix},
\end{equation*}
where $I_k$ is the identity matrix $k\times k$ and $L$ is the Cholesky factor of the sparse matrix $S\Lambda_S S^T +FF^T$. 

Let $G:=L^{-1}D\Lambda_D^{1/2}$ and $J:=-L^{-1}F$. From the definition of $F$, it follows that $J$ is a matrix produced by $d$ columns of $L^{-1}$ times a constant, making it a parallelizable operation from the perspective of the computational implementation. 

The preconditioned matrix is
\begin{equation*}
     M^{-1}KM^{-T}=\begin{bmatrix}
    I_m+GG^{T}-JJ^T&&\\
     & -\Lambda_D^{-1}&\\
     &&I_d
    \end{bmatrix}.
\end{equation*}

To recover the solution, that is, the search Newton's direction $\Delta y$ defined in \eqref{eqnormal}, we must follow some steps summarized in the Algorithm \ref{newprecond}:

\begin{algorithm}[H]
\begin{algorithmic}[1]
     \STATE \textbf{Input:} $A=[S,D], \Lambda = [\Lambda_S, \Lambda_D]$, $q$.
     \STATE \textbf{Output: } Newton's direction $\Delta y$
    \STATE Compute the Cholesky factorization $LL^T=S\Lambda_S S^T+FF^T$
    \STATE Solve $LG = D\Lambda_D ^{1/2}$, $G:m\times k$. Solving $k$ triangular linear systems.
    \STATE Solve $LJ = -F$, $J: m\times d$. Solving $d$ triangular linear systems.
    \STATE Solve the triangular system $L\hat{q}_1 = q$
    \STATE Solve with iterative methods the linear system $(I_m+ GG^T -JJ^T)\hat{\omega}_1=\hat{q}_1$
    \STATE Solve the triangular system $L^T \Delta y=\hat{\omega}_1$.
    \end{algorithmic}
    \caption{Algorithm to get the search direction $\Delta y$}
    \label{newprecond}
\end{algorithm}

\subsection{Numerical stability}

In this section, we provide a detailed explanation of asymptotic notation and its application to matrices. Suppose $A$ is well-conditioned and well-scaled. For asymptotic properties related to each type of boundary, we refer to (\cite{algorithms}). Let $B$ be a matrix dimension $m \times n$,  such that depends on a variable $\gamma$. We introduce some definitions of asymptotic behavior. 

\begin{definition}
We say that $B$ is \textit{uniformly bounded} if there exists a positive constant $c$ that is not dependent on $\gamma$, such that the absolute value of each entry of the matrix $B$ is upper bounded by $c$. We denote this as $B=\mathcal{O}(1)$. This means the matrix $B$ remains relatively constant as $\gamma$ changes.
\end{definition}

\begin{definition}
We say $B=\mathcal{O}(\gamma)$ if there is a positive constants $c$ such that each entry of $B$ satisfies $|B_{ij}| \leq c  \gamma$.
\end{definition}

\begin{definition}
    We say $B=\Theta(\gamma)$ if there are two positive constants $c_1$ and $c_2$, such that each entry of $B$ satisfies $c_1 \gamma \leq |B_{ij}| \leq c_2  \gamma$. This means that the entries of the matrix $B$ are growing at the same rate as $\gamma$. Similar definitions are used for the notation $B=\Theta(1)$.
\end{definition}

It follows from definition that if a matrix $B=\Theta(\gamma)$, then, $B=\mathcal{O}(\gamma)$.

\begin{remark}
    From here, we analyze the case where a matrix $M$ depends on $\mu=x^Tz+v^Tw$. In other words, we examine the uniformly boundless case in terms of the gap $\mu$.
\end{remark}

\begin{remark}
\label{remarkgamma}
    Moving on to the optimality conditions, we note that the gap $\mu=x^Tz+v^Tw$ is approaching zero as the Interior-Point Method is converging to an optimal solution. Since we are working on primal-dual methods, then, the diagonal entries of scale diagonal matrix $\Lambda$ are either $\Theta(\mu)$ or $\Theta(\mu^{-1})$ (See \cite{wright1997primal}). 
\end{remark}

Suppose that $F$ is $\mathcal{O}(1)$.

The diagonal matrix $\Lambda= Z^{-1}X+W^{-1}V$ depends on the values of primal and dual variables, i.e., it depends on $\mu$. We are interested in the numerical stability of the matrix $W=I_m+GG^{T}-JJ^T$. Specifically, we want to determine if $W$ is uniformly bounded as the Interior-Point Method approaches an optimal solution, as indicated by the convergence of $\mu$ to zero.

 The next two theorems from D. Goldfarb and  K. Scheinberg in (\cite{Choleskyproduct}) are important to analyze the numerical stability of matrix $W$. 

\begin{theorem}[\cite{Choleskyproduct}]
	\label{Lisunifbounded}
	Let $A:m\times n$ and $\Lambda$ be a positive diagonal matrix, whose entries depend on $\mu$ and $\tilde{L}\Sigma \tilde{L}^{T}$ be the Cholesky factorization of $A\Lambda A^{T}$, where $\tilde{L}$ is a lower triangular matrix with ones on the main diagonal and $\Sigma$ is a positive diagonal matrix. Then, the entries of $\tilde{L}$ are uniformly bounded when $\mu$ is approaching zero.	
\end{theorem}

\begin{theorem}[\cite{Choleskyproduct}]
	\label{inverseunifbounded}
	If $\tilde{L}$ is a lower triangular matrix with ones on the main diagonal, and the subdiagonal entries that depend on a parameter $\mu$ are such that they are uniformly bounded, then the entries of inverse $\tilde{L}^{-1}$ are also uniformly bounded. 
\end{theorem}

 To analyze the asymptotic behavior of $W$, we proceed to analyze the behavior of $LL^T=S\Lambda_S S^T+FF^T$ in the following proposition.

\begin{proposition}
\label{STS}
    $LL^T= \mathcal{O}(\gamma^{-1})$.
\end{proposition}
\begin{proof}
    Let $i,j\in \{ 1,\dots ,m \}$, then,
    \begin{equation*}
        (S\Lambda_S S^T)_{ij} = \sum_{l=1}^{n-k} \Lambda_{S_l} S_{il}S_{jl}.
    \end{equation*}
    According to Remark \ref{remarkgamma}, the entries of scale matrix $\Lambda_S$ are either $\Theta(\mu)$ or $\Theta(\mu^{-1})$, therefore, the nonzero entries of $S\Lambda_S S^T$ are a linear combination of order $\Theta(\mu)$ or $\Theta(\mu^{-1})$.  
    
    Since we are analyzing the case of $\mu\to 0^+$ and there is at least one entry in the matrix $S\Lambda_S S^T$ with order $\Theta(\mu^{-1})$, therefore, by asymptotic properties, we obtain that  $S\Lambda_S S^T$ is $\mathcal{O}(\mu^{-1})$.

    Since we supposed that the nonzero entries of $F$ are $\mathcal{O}(1)$ and $\mu\to 0^+$, thus, by asymptotic properties, the entries of $S\Lambda_S S^T +FF^T$ are $\mathcal{O}(\mu^{-1})$. Therefore, $LL^T= \mathcal{O}(\mu^{-1})$.
    
\end{proof}

The previous Proposition \ref{STS} also confirms the known fact that normal equations are ill-conditioned when the Primal-Dual method is approaching an optimal solution (\cite{wright1997primal}).

Once we have analyzed the behavior of $LL^T$, we present the following results:

\begin{lemma}
\label{lemmaa}
    Let $LL^T=S\Lambda_S S^T+FF^T$. There is a lower triangular matrix $\tilde{L}$ with ones in the main diagonal and $\Sigma$ positive diagonal matrix such that $L=\tilde{L}\Sigma^{1/2}$. In particular, the diagonal entries of $\Sigma$ are $\Theta(\mu^{-1})$. Thus, $\Sigma=\mathcal{O}(\mu^{-1})$.
    
\end{lemma}
\begin{proof}
    It is trivial that we obtain the factorization $L=\tilde{L}\Sigma^{1/2}$. From Proposition \ref{STS}, $LL^T=\mathcal{O}(\mu^{-1})$. Then, $\tilde{L}\Sigma\tilde{L}^T=\mathcal{O}(\mu^{-1})$. From Theorem \ref{Lisunifbounded} we obtain that $\tilde{L}$ and $\tilde{L}^T$ are uniformly bounded matrices. Therefore, the positive diagonal matrix is $\Sigma=\mathcal{O}(\mu^{-1})$. 
    Since we are working on Primal-Dual method, from Proposition \ref{STS} and Remark \ref{remarkgamma}, it follows that the matrix $\Sigma$ has diagonal entries of order $\Theta(\mu)$ and $\Theta(\mu^{-1})$. Since we are analyzing the case $\mu\to 0^+$, we conclude that the main diagonal of $\Sigma$ is $\Theta(\mu^{-1})$.
\end{proof}


\begin{proposition}
\label{GGTprop}
    $GG^T=\mathcal{O}(1)$.
\end{proposition}
\begin{proof}
In a similar way to Proposition \ref{STS}, for $i,j\in \{ 1,\dots, m \}$ we have 
 \begin{equation*}
        (D\Lambda_D D^T)_{ij} = \sum_{l=1}^{k} \Lambda_{D_l} D_{il}D_{jl}.
    \end{equation*}

It implies that matrix $D\Lambda_D D^T$ is $\mathcal{O}(\mu^ {-1})$.

On the other hand, from Lemma \ref{lemmaa}, $S\Lambda_S S^T+FF^T = \tilde{L}\Sigma \tilde{L}^T$, and diagonal matrix $\Sigma$ is $\mathcal{O}(\mu ^{-1})$, with positive diagonal entries of order $\Theta(\mu^{-1})$. Also, from Theorem \ref{inverseunifbounded}, we obtain $\tilde{L}^{-1}$ is uniformly bounded.

It implies that the diagonal entries of $\Sigma^{-1}$ are of order $\Theta(\mu)$. Thus, there are positive constants $c_1,c_2$ such that $c_1 \mu\leq \Sigma_{ii}^{-1}\leq c_2 \mu$. It follows that
\begin{equation*}
    \sqrt{c_1}\mu^{1/2} \leq \frac{1}{\sqrt{\Sigma_{ii}}} \leq \sqrt{c_2}\mu^{1/2}.
\end{equation*}
Then, the diagonal entries of $\Sigma^{-1/2}$ are $\Theta(\mu^{1/2})$.

Therefore, combining those results and using the facts that $D\Lambda_D D^T=\mathcal{O}(\mu^{-1})$, then, for $GG^T$ we obtain
\begin{equation*}
       \underbrace{\Sigma^{-1/2}}_{\mathcal{O}(\mu^{1/2})} \underbrace{\tilde{L}^{-1}}_{\mathcal{O}(1)} \underbrace{D \Lambda_D D^T}_{\mathcal{O}(\mu^{-1})} \underbrace{\tilde{L}^{-T}}_{\mathcal{O}(1)} \underbrace{\Sigma^{-1/2}}_{\mathcal{O}(\mu^{1/2})} = GG^T.
   \end{equation*}

This is, by asymptotic properties, $GG^T=\mathcal{O}(1)$.
     
\end{proof}

Similarly, we may analyze the matrix $JJ^T$.

\begin{proposition}
\label{JJTprop}
    $JJ^T=\mathcal{O}(1)$.
\end{proposition}
\begin{proof}
    According to Lemma \ref{lemmaa}, we get that the inverse of Cholesky factor $L^{-1}=\Sigma^{-1/2}\tilde{L}^{-1}$ is $\mathcal{O}(\mu^{1/2})$. Since the entries of $F$ have order $\mathcal{O}(1)$, then, $J= L^{-1}F= \mathcal{O}(\mu^{1/2})$. Thus, $JJ^T=\mathcal{O}(\mu)$. Given that we are analyzing the case $\mu\to 0^+$, we conclude that $JJ^T$ is uniformly bounded,  this is $JJ^T=\mathcal{O}(1)$.
\end{proof}

   We finally conclude about the asymptotic behavior of $W$ when Interior-Point Methods is approaching an optimal solution.

   \begin{theorem}
    $W$ is uniformly bounded.
\end{theorem}
\begin{proof}
    Let $i,j\in \{ 1,\dots, m\}$, then 
    \begin{equation*}
		(W)_{ij}=
	\begin{cases}
		(GG^T)_{ij}-(JJ^T)_{ij}&\text{ if } i\neq j,\\
		1+(GG^T)_{ij}- (JJ^T)_{ij} & \text{ if } i= j.
	\end{cases}
\end{equation*}
    Considering that matrices $GG^T$ and $JJ^T$ are uniformly bounded (Propositions \ref{GGTprop} and \ref{JJTprop}), there are positive constants $\alpha, \beta$ for all $i,j \in \{ 1,\dots, m \}$  such that $|(JJ^T)_{ij}|\leq \alpha$ and $|(GG^T)_{ij}|\leq \beta$. Therefore for each entry of $W$ we obtain   
\begin{align*}
    |W_{ij}|&=|(GG^T)_{ij}-(JJ^T)_{ij}|\leq \alpha +\beta & \qquad \mbox{if } i\neq j,\\
    |W_{ij}|&=|1+(GG^T)_{ij}-(JJ^T)_{ij}|\leq \alpha + \beta +1  \qquad &\mbox{if } i = j.
\end{align*}

Anyway, $|W_{ij}|\leq \alpha + \beta +1 $. Thus, $W$ is uniformly bounded when Interior-Point Methods is approaching an optimal solution. 

\end{proof}

\begin{remark}
    Since we suppose that the constraint matrix is well-conditioned and well-scaled, the positive constants associated with the Sigma ($\Sigma(\cdot)$) and Big-O ($\mathcal{O}(\cdot)$) notation are moderate (See \cite{Choleskyproduct}).
\end{remark}

\subsection{Properties of the preconditioned system}
\label{propertiesofW}

The present section is focused on the theoretical properties of matrix $W=I_m+GG^{T}-JJ^T$. 

\begin{proposition}
\label{eigenvalue1}
	$W=I_m+GG^{T}-JJ^T$ is positive definite.
 \end{proposition}
 \begin{proof}
     Based on the definition of $J$ and $G$ we obtain for $W$:
\begin{align}
\label{precondisA}
W &= I+GG^T-JJ^T \nonumber
\\ &= I+L^{-1}D\Lambda_DD^TL^{-T}-L^{-1}FF^TL^{-T} \nonumber\\
& = L^{-1}(LL^T +D\Lambda_D D^T -FF^T)L^{-T}\nonumber\\
& = L^{-1} (A\Lambda A^T)L^{-T}.
\end{align}

Since $A$ is full-row rank, $\Lambda$ is a positive diagonal matrix and $L$ is the Cholesky factor of the sparse matrix $S\Lambda_S S^T +FF^T$, therefore $W$ is a positive definite matrix.

 \end{proof}

\begin{proposition}
	\label{eigmax}
The largest eigenvalue of $W$ satisfies	$\lambda_{\max} (W) = \| L^{-1}A\Lambda^{1/2}\|_2^2$.
\end{proposition}
\begin{proof}
	According to the definition $G=L^{-1}D\Lambda_D^{1/2}$ and $J=-L^{-1}F$, then we have
	\begin{equation*}
		\begin{split}
            \lambda_{\max}(W)&=
			\lambda_{\max}(I+GG^T-JJ^T)\\
            &=\lambda_{\max}(I+(L^{-1}D\Lambda_D^{1/2})(L^{-1}D\Lambda_D^{1/2})^{T}-(L^{-1}F)(L^{-1}F)^T)\\& 
			= \lambda_{\max}(I+L^{-1}D\Lambda_DD^{T}L^{-T}-L^{-1}FF^{T}L^{-T})\\
			&=\lambda_{\max}(L^{-1}(LL^{T}+D\Lambda_DD^{T}-FF^T)L^{-T})\\
			&=\lambda_{\max}(L^{-1}(A\Lambda A^{T})L^{-T})\\
			&=\lambda_{\max}(L^{-1}(A\Lambda^{1/2}\Lambda^{1/2} A^{T})L^{-T})\\
			&=\lambda_{\max}((L^{-1}A\Lambda^{1/2})(L^{-1}A\Lambda^{1/2})^{T} ).
		\end{split}	
	\end{equation*}
Based on the $2$-norm equivalences, we got
\begin{equation*}
	\begin{split}
	\lambda_{\max} (W)=	\lambda_{\max}((L^{-1}A\Lambda^{1/2})(L^{-1}A\Lambda^{1/2})^{T} )= \| L^{-1}A\Lambda^{1/2}\|_2^2.
	\end{split}
\end{equation*}
\end{proof}

The Proposition \ref{eigmax} states an equivalence of $\lambda_{\max}$ that depends on diagonal matrix $\Lambda$. Thus, we proceed to obtain an upper boundary of $\lambda_{\max}$ independent of $\mu$.

\begin{theorem}
\label{largesteigenvalue}
    There is a constant $c>0$ such that $\lambda_{\max}(W)\leq c$.
\end{theorem}

\begin{proof}
The result follows from the fact that $W$ is uniformly bounded. There is a positive constant $c_W$ such that for each $i,j\in \{ 1,\dots , m\}$,  $|W_{ij}|\leq c_W$. Then, if $(\lambda_{\max},v)$ is an eigenpair, such that $v$ is a unit vector, then
\begin{align*}
    \lambda_{\max} = |\lambda_{\max}|=|\langle v, Wv \rangle| =\left| \sum_{i=1}^m v_i(Wv)_i \right| ,
\end{align*}
where 
\begin{equation*}
    (Wv)_i = \sum_{j=1}^{m} W_{ij}v_j.
\end{equation*}

Therefore, applying triangular inequality 
\begin{align*}
    \lambda_{\max}& \leq \sum_{i=1}^m |v_i| \left|\sum_{j=1}^{m} W_{ij}v_j\right|
    \leq \sum_{i=1}^m |v_i| \sum_{j=1}^{m} |W_{ij}||v_j|
    \leq \sum_{i=1}^m |v_i| mc_W \leq c_Wm^2.
\end{align*}
\end{proof}

Since $W$ is a normal matrix, thus, the condition number $\kappa(W)$ satisfies (See \cite{meyer})
\begin{equation*}
    \kappa(W) = \frac{\lambda_{\max}(W)}{\lambda_{\min}(W)},
\end{equation*}
where $\lambda_{\min}(W)$ is the minimum eigenvalue of W.

Given that $W$ is a positive definite matrix, there is $\epsilon >0$ such that $\lambda_{\min}\geq \epsilon$. Therefore, from Theorem \ref{largesteigenvalue} we get an upper bound for the condition number
\begin{equation*}
    \kappa(W) \leq \frac{m^2c_W}{\epsilon}.
\end{equation*}

\begin{proposition}
    If $F=0$, then $\lambda_{\min} (W) \geq 1$.
\end{proposition}
\begin{proof}
    Let $(\lambda,v)$ be an eigenpair, such that $v$ is unit vector. Then, if $F=0$
    \begin{align*}
        \lambda v=(I+GG^T)v=v+GG^Tv.
    \end{align*}
Thus, it follows
\begin{align*}
    \lambda &= \lambda \langle v,v \rangle = \langle v,(I+GG^T)v \rangle
    = \langle v,v \rangle +\langle v,GG^Tv \rangle = 1+\| G^Tv\|_2^2.
\end{align*}
Since $\| G^Tv \|\geq 0$, then $\lambda \geq   1$. It is valid for any eigenvalue, in
particular $\lambda_{\min}\geq 1$.

\end{proof}

From the previous proposition, we conclude that if $S$ is full-row rank, this is, $F=0$, 
\begin{equation*}
    \kappa (W) \leq \lambda_{\max}(W).
\end{equation*}

\section{Computational experiments}
\label{experim}

In this section, we present computational experiments for linear programming problems involving dense columns that can be addressed using the proposed preconditioner. Before exposing the computational experiment results, we present some important facts about the implementation.

\subsection{PCx}

PCx is an Interior-Point predictor-corrector linear programming package whose theoretical results are in (\cite{PCx}). The Interior-Point is developed in C and the Cholesky factorization in Fortran 77, using the approach of sparse factorization of (\cite{NgPeyton, NgPeytonoriginal}).

\subsection{Implementation}
\label{implem}
Through a modification of the PCx code in (\cite{PCx}), the proposed preconditioner for handling linear programming problems with dense columns in constraint matrices is developed in C and Fortran 77. The modification focuses on search direction $\Delta y$ (output in Algorithm \ref{newprecond}) and the Cholesky factorization, to include the matrix $F$.

\subsection{Permutations matrices}

Two permutation matrices, $P_1$ and $P_2$, are required for the preprocessing procedure applied to the constraint matrix. The first one seeks to split the constraint matrix into sparse and dense parts, i.e., $AP_1=[S, D]$ and $P_1 \Lambda P_1^T = \textit{diag}(\Lambda_S, \Lambda_D)$.

The row permutation matrix is associated with the Cholesky factorization. In the present implementation, we use the process developed by Ng and Peyton in (\cite{NgPeyton, NgPeytonoriginal}), which focuses on factorization for sparse matrices and employs the Minimum Degree Heuristics (\cite{mindegree}). In the proposed preconditioner, our goal is to find the Cholesky factorization 
\begin{equation*}
    P_2(S\Lambda_S S^T)P_2^T =LL^T.
\end{equation*}

It is clear that even if the constraint matrix is preprocessed, it may contain some redundant rows. In addition, one of the columns associated with the dense matrix $D$ may be a possible pivot to obtain the Cholesky factorization of the normal matrix $A\Lambda A^T$; resulting in small pivots during the computation of the Cholesky factor of $S\Lambda_S S^T$. 

For the present implementation, we consider a small pivot if 
\begin{equation*}
    M_{ii}^{(i-1)}\leq 10^{-30} \max_{j=1,\dots, m} (S\Lambda_S S^T)_{jj}^2,
\end{equation*}
where $M^{(i-1)}$ is the submatrix after $(i-1)$ iterations of the Cholesky process. 

At each small pivot position, we add up a positive constant in the pivot column where the matrix factorization has failed. For this, a matrix $F$ with $d\leq k$ columns is formed by canonical vectors, each one associated with a small pivot. Finally, we obtain the Cholesky factorization
\begin{equation}
\label{cholperm}
    P_2 S\Lambda_S S^T P_2^T +FF^T= LL^T.
\end{equation}

Therefore, after applying the Minimum Degree heuristic, the permuted normal equation matrix is
\begin{align*}
    P_2 A \Lambda A^T P_2^T &= P_2 S \Lambda_S S^T P_2^T + FF^T + P_2 D \Lambda_D D^T P_2^T -FF^T \nonumber\\
    & = LL^T + P_2 D \Lambda_D D^T P_2^T -FF^T.
\end{align*}

The process to compute the search direction $\Delta y$ of the new preconditioner is in Algorithm \ref{newprecondperm} as follows

\begin{algorithm}[H]
    \begin{algorithmic}[1]
    \STATE \textbf{Input:} $AP_1=[S,D], P_1\Lambda P_1 = [\Lambda_S, \Lambda_D]$, $q$.
    \STATE \textbf{Output: } Direction $\Delta y$
        \STATE Compute the Cholesky factorization $LL^T=P_2S\Lambda_S S^TP_2^T+FF^T$
        \STATE Solve $LG = P_2D\Lambda_D ^{1/2}$, $G:m\times k$
        \STATE Solve $LJ = -F$, $J: m\times d$
        \STATE Solve the triangular system $L\hat{q}_1 = P_2q$
        \STATE Solve with iterative methods the linear system $(I_m+ GG^T -JJ^T)\hat{\omega}_1=\hat{q}_1$
        \STATE Solve the triangular system $L^T\tilde{\Delta }y=\hat{w}_1$
        \STATE Compute permutation $P_2\Delta y = \tilde{\Delta}y$
    \end{algorithmic}
    \caption{Algorithm to get the search direction $\Delta y$ with permutations}
    \label{newprecondperm}
\end{algorithm}

\subsection{Dense criteria}
\label{densecrit}

For the present preconditioner, we use the dense criteria determined by PCx in (\cite{PCx}), that is a more rigorous criterion in comparison to the Andersen proposed in (\cite{andersen}). For more detail, see \cite{PCx} or website \url{https://pages.cs.wisc.edu/~swright/PCx/}.
\newline

Once the implementation requirements are defined, we proceed to expose the computational results. The computational experiments are performed in a Linux environment on an Intel Core i7-3770K processor running at $3.50$GHz, with $32$ GB RAM. The numerical experiments apply the preconditioner defined in Section \ref{sec3} and a modified version of the PCx algorithm (\cite{PCx}). The Interior-Point method is implemented in C, and the Cholesky factorization associated with the sparse matrix $S\Lambda_S S^T +FF^T$ is developed in Fortran 77, following the approach in (\cite{NgPeyton}). For presolving process, we employed the Gurobi Optimization Software (\cite{gurobi2020gurobi}) and the PCx (\cite{PCx}).
\newline

Our proposal, which we denote as PCx\_mod, is compared with the original PCx package, in which we denote PCx\_orig. Also, we compare it with the Modified Schur-Complement of Andersen (PCx\_Andersen) in \cite{andersen}. We use the same density criteria for the three cases for a fair comparison. Likewise, we also analyze the case of not splitting dense and sparse columns (PCx\_full). All of these are developed with a modification of PCx.

In the case of PCx\_orig, the search direction $\Delta y$ for the Interior-Point methods is computed with the Sherman-Morrison-Woodbury formula. If this solution does not meet a certain tolerance, the Preconditioned Conjugate Gradient method is applied. This method uses the Cholesky factor $L$ of the sparse component as a preconditioner. It is noteworthy that the PCx approach solves the normal equations \eqref{eqnormal} for each predictor and corrector step, using the corresponding right-hand side. For the tested problems, the Preconditioned Conjugate Gradient method was not necessary.

Similarly, for PCx\_full, the search direction $\Delta y$ is computed using the same method as PCx\_orig, but the density criterion is modified to ensure that there are always zero dense columns.

On the other hand, for PCx\_Andersen, the search direction $\Delta y$  is computed by applying LU factorization as proposed in \cite{andersen}.
\newline 

For the present paper, most of the tested Linear Programming problems are from the paper (\cite{castro2006minimum}), for the particular case of Minimum-distance controlled perturbation methods with L-infinity norm (library in \url{http://www-eio.upc.es/~jcastro/}).

The $22$ cases that are tested and their descriptions are summarized in Table \ref{datosproblemas}. $m$ and $n$ represent the number of rows and columns respectively of constraint matrix $A$, Nnz(A) represents the nonzero entries of $A$ and ndense is the number of dense columns of $A$. It is highlighted that the table only includes the size of each Linear Programming (LP) problem after applying the presolving process. For all the examples of this table we get $F=0$, that is, the matrix $S\Lambda_S S^T$ is full rank.

\begin{table}[H]
\caption{Tested problems of J. Castro database.}
\label{datosproblemas}
\begin{tabular*}{\hsize}{@{}@{\extracolsep{\fill}}lrrrr@{}}
\hline
  LP problem & {$m$}  & {$n$}  & {Nnz(A)}  & {ndense}   \\
\hline
Linf\_bts4          & 55145                 & 63808                 & 265131                     & 1                          \\
Linf\_cbs           & 3184                  & 3354                  & 12652                      & 1                          \\
Linf\_five20c       & 47128                 & 48949                 & 233628                     & 1                          \\
Linf\_five20b       & 48143     &50690                 & 239935                     & 1                          \\
Linf\_hier13        & 3176                  & 3265                  & 14859                      & 1                          \\
Linf\_hier13x7x7d   & 722                   & 825                  & 3404                       & 1                          \\
Linf\_hier13x13x7d  & 1472                  & 1574                  & 6621                       & 1                          \\
Linf\_hier13x13x13a & 2992                  & 3096                  & 13293                      & 1                          \\
Linf\_hier13x13x13b & 2992                  & 3096                  & 13293                       & 1                          \\
Linf\_hier13x13x13c & 2992                  & 3096                  & 13293                      & 1                          \\
Linf\_hier13x13x13d & 2992                  & 3096                  & 13293                       & 1                          \\
Linf\_hier13x13x13e & 2992                  & 3096                  & 13293                      & 1                          \\
Linf\_jjtabeltest3  & 2563                  & 2554                  & 9770                      & 2                          \\
Linf\_nine5d        & 9853                  & 10532                 & 47308                      & 1                          \\
Linf\_nine12        & 15815                 & 17592                 & 80293                      & 1                          \\
Linf\_ninenew       & 9185                  & 10072                 & 46443                      & 1                          \\
Linf\_osorio        & 10199                 & 19995                 & 69779                      & 1                          \\
Linf\_table1        & 1308                  & 1775                  & 7197                       & 1                          \\
Linf\_table3        & 3477                  & 4165                  & 18843                      & 1                          \\
Linf\_table6        & 1308                  & 1775                  & 7197                       & 1                          \\
Linf\_table8        & 1342                  & 2539                  & 8800                      & 1                          \\
Linf\_two5in6  & 5612                  & 5989                  & 27129                      & 1                        \\
\hline
\end{tabular*}
\end{table}

Table \ref{iterations} summarizes the total iterations to obtain an optimal solution for each linear programming problem. 
The ``UNKNOWN" result is due to the convergence rate being very slow or the $\mu$ value increasing from one iteration to another; ``-" refers to a problem takes more than 5 hours to be solved. 

While PCx\_full may appear to require fewer iterations in approximately $45\%$ of cases, a subsequent table reveals that the time required is considerably higher, as anticipated. Conversely, the results for PCx\_mod, PCx\_orig, and PCx\_Andersen are comparable, obtaining $41\%$ of LP problems with the same required iterations. In $18\%$ of tested problems, our proposed method exhibits fewer iterations. However, the robustness of our approach becomes evident as it demonstrates the capability to solve the problem \textit{Linf\_table3}, marked as unknown or infeasible to the other approaches. PCx\_Andersen failed in $23\%$ of problems.

\begin{table}[h]
\caption{Total iterations for each approaching method}%
\label{iterations}
\begin{tabular*}{\textwidth}{@{}@{\extracolsep{\fill}}lrrrrr@{}}
\hline
LP problem & Gond C. & \multicolumn{4}{c}{Number of iterations}  \\[-6pt]
&  & \multicolumn{4}{l}{\hrulefill}  \\
 &  & PCx\_mod & PCx\_orig & PCx\_Andersen & PCx\_full \\\hline
Linf\_bts4         & 3                             & 39                   & 40                             & \textbf{34}                                 & 36                             \\
Linf\_cbs          & 2                             & 18                            & 18                             & 18                                 & 21                             \\
Linf\_five20c      & 10                            & 20                            & 20                             & 20                                 & 20                             \\
Linf\_five20b      & 6                             & \textbf{20}                   & \textbf{20}                             & 21                                 & -                              \\
Linf\_hier13        & 1                             & \textbf{13}                   & 14                             & \textbf{13}                        & 15                             \\
Linf\_hier13x7x7d   & 0                             & 19                            & 19                             & 19                                 & 19                             \\
Linf\_hier13x13x7d  & 0                             & 20                            & 20                             & UNKNOWN        & \textbf{17}                             \\
Linf\_hier13x13x13a & 1                             & 20                            & 20                             & 20                                 & \textbf{18}                             \\
Linf\_hier13x13x13b & 1                             & 20                            & 20                             & 20                                 & \textbf{18}                             \\
Linf\_hier13x13x13c & 1                             & 20                            & 20                             & 20                                 & \textbf{18}                             \\
Linf\_hier13x13x13d & 1                             & 20                            & 20                             & UNKNOWN        & \textbf{18}                             \\
Linf\_hier13x13x13e & 1                             & 20                            & 20                             & UNKNOWN       & \textbf{18}                             \\
Linf\_jjtabeltest3  & 0                             & 28                            & 28                             & INFEASIBLE     & \textbf{24}                             \\
Linf\_nine5d        & 2                             & 16                            & 16                             & 16                                 & \textbf{15}                             \\
Linf\_nine12        & 4                             & 14                            & 14                             & 14                                 & 15                             \\
Linf\_ninenew       & 4                             & \textbf{13}                   & 14                             & 14                                 & 14                             \\
Linf\_osorio       & 0                             & 16                   & 17                             & 16                        & 15                             \\
Linf\_table1        & 1                             & 26                            & 25                             & \textbf{23}                        & 25                             \\
Linf\_table3        & 2                             & 27                   & INFEASIBLE                             & UNKNOWN        & \textbf{24}                             \\
Linf\_table6        & 0                             & \textbf{22}                   & 28                             & \textbf{22}                        & 25                             \\
Linf\_table8       & 0                             & 18                            & 18                             & 18                                 & 19                             \\
Linf\_two5in6  & 1                             & 21                   & 22                             & 21                        & \textbf{20}\\
\hline
\end{tabular*}
\end{table}
\hspace{5mm}

The running time for each case is summarized in Table \ref{runningtime}. The values in the table do not include the time for the preprocessing process. We can observe that $45\%$ problems produce results that are essentially the same for PCx\_mod, PCx\_orig, and PCx\_Andersen. Our proposal had better success in $37\%$ cases. The case of not considering sparse and dense columns separated, denoted as PCx\_full required extremely high times in comparison to other methods as anticipated, requiring more than $10$ times the running time for all tested problems.

\begin{table}[H]
\caption{Running time of each approaching method}%
\label{runningtime}
\begin{tabular*}{\textwidth}{@{}@{\extracolsep{\fill}}lrrrr@{}}
\hline
LP problem  & \multicolumn{4}{c}{Time(s)}  \\[-6pt]
&   \multicolumn{4}{l}{\hrulefill}  \\
   & PCx\_mod & PCx\_orig & PCx\_Andersen & PCx\_full \\\hline
 Linf\_bts4                       & \textbf{23.01}                & 25.86                          & 24.26                     & 110446.02                      \\
Linf\_cbs                        & \textbf{0.43}                          & \textbf{0.43}                           & 0.44                               & 49.5                           \\
Linf\_five20c                    & \textbf{1064.38}              & 1064.44                        & 1065.45                            & 51763.45                       \\
Linf\_five20b                    & \textbf{189.13}               & 196.47               & 189.95                             & -           \\
Linf\_hier13                     & \textbf{0.14}                 & \textbf{0.14}                           & \textbf{0.14}                               & 10.62                          \\
Linf\_hier13x7x7d                & 0.02                          & 0.02                           & 0.02                               & 0.24                           \\
Linf\_hier13x13x7d               & \textbf{0.03}                 & \textbf{0.03}                  & UNKNOWN       & 1.18                           \\
Linf\_hier13x13x13a              & 0.16                          & 0.16                           & 0.16                               & 11.88                          \\
Linf\_hier13x13x13b              & 0.16                          & 0.16                           & 0.16                               & 11.83                          \\
Linf\_hier13x13x13c              & 0.15                          & 0.15                           & 0.16                               & 11.96                          \\
Linf\_hier13x13x13d              & 0.16                          & 0.16                           & UNKNOWN       & 11.9                           \\
Linf\_hier13x13x13e              & 0.16                          & 0.16                           & UNKNOWN       & 12.7                           \\
Linf\_jjtabeltest3               & 0.06                          & \textbf{0.05}                           & INFEASIBLE     & 7.33                           \\
Linf\_nine5d                     & \textbf{1.76}                 & 1.77                           & 1.77                               & 334.5                          \\
Linf\_nine12                     & \textbf{17.44}                & 17.69                          & 17.91                              & 1769.12                        \\
Linf\_ninenew                    & \textbf{13.58}                & 13.91                          & 14.24                              & 406.35                         \\
Linf\_osorio                     & 0.15                 & 0.15                           & 0.15                      & 1561.2                         \\
Linf\_table1                     & 0.06                          & 0.05                           & \textbf{0.04}                      & 1.8                            \\
Linf\_table3                     & \textbf{1.11}                 & INFEASIBLE & UNKNOWN       & 34.91                          \\
Linf\_table6                     & 0.05                          & 0.05                           & 0.05                               & 1.8                            \\
Linf\_table8                     & 0.02                          & 0.02                           & 0.02                               & 4.52                           \\
Linf\_two5in6               & \textbf{0.58}                 & 0.6                            & \textbf{0.58}                      & 89.4    
\\
\hline
\end{tabular*}
\end{table}
\hspace{4mm}

\begin{figure}[h]
\centering
\includegraphics[width=.65\textwidth]{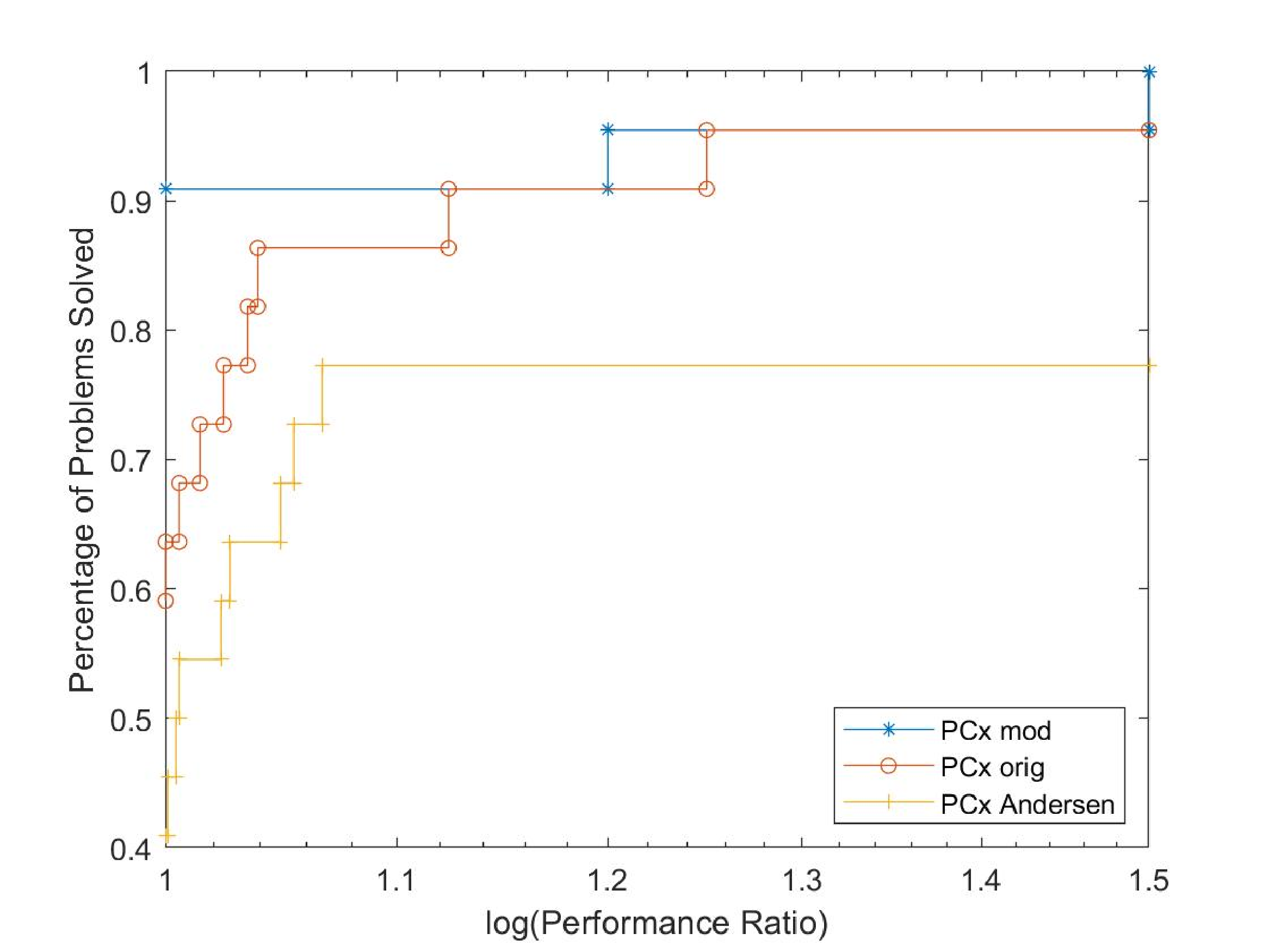}
\caption{ {\small Performance profile for the three dense approaches.}}
\label{figura02}
\end{figure}

Figure \ref{figura02} shows the performance profile illustrating the time taken to solve 22 problems across three different approaches for solving LP problems. Specifically, it compares the performance of these approaches: PCx\_mod, PCx\_orig, and PCx\_Andersen according to Tables \ref{runningtime}. The x-axis represents the performance ratio in the log base, while the y-axis shows the probability that each solver achieves a certain performance level relative to the best solver across all problems.

We observe that there is a $90\%$ probability to be the fastest method. This observation underscores the efficiency and effectiveness of PCx\_mod in solving the given LP problems with dense columns.

In Table \ref{cgproblemas}, we present the total iterations required to solve the linear system associated with the matrix $W=I_m+GG^T-JJ^T$ using Conjugate Gradient. The columns ``Predictor step" and ``Corrector step" represent the mean number of iterations in the Conjugate Gradient Iterations method. For the decimal values, this implies that, for certain iterations of the Interior-Point Method (IPM), required more than one Conjugate Gradient (CG) iteration. The column ``Total CG iterations" represents the cumulative number of iterations needed to obtain an optimal solution, adding all iterations for each step of the IPM.

The values in Table \ref{cgproblemas} involve the required iterations to compute the search direction $\Delta y$ according to step $7$ in Algorithm \ref{newprecond}. It is observed that approximately two Conjugate Gradient iterations are required to solve each iteration of the Interior Point Method, one for the predictor step and the other for the corrector. It is noteworthy that the other solutions apply direct techniques for solving the linear system. Additionally, we emphasize that, due to the small number of iterations required, preconditioners for iterative methods are not necessary. 

\begin{table}[h]
\caption{Conjugate gradients results of PCx\_mod}\label{cgproblemas}%
\begin{tabular*}{\textwidth}{@{}@{\extracolsep{\fill}}lrrr@{}}\hline
LP Problem & Total CG iterations & \multicolumn{2}{c}{Mean of CG iterations}  \\
&&\multicolumn{2}{c}{\hrulefill} \\
 & & Predictor step & Corrector step  \\\hline
Linf\_bts4                       & 76                                                        & 1                                   & 1                                   \\
Linf\_cbs                        & 34                                                        & 1                                   & 1                                   \\
Linf\_five20c                    & 38                                                        & 1                                   & 1                                   \\
Linf\_five20b                    & 38                                                        & 1                                   & 1                                   \\
Linf\_hier13                     & 24                                                        & 1                                   & 1                                   \\
Linf\_hier13x7x7d                & 36                                                        & 1.11                                   & 1.11                                  \\
Linf\_hier13x13x7d               & 40                                                        & 1                                   & 1                                   \\
Linf\_hier13x13x13a              & 39                                                        & 1.05                                   & 1                                   \\
Linf\_hier13x13x13b             & 39                                                        & 1.05                                   & 1                                   \\
Linf\_hier13x13x13c              & 39                                                        & 1.05                                   & 1                                  \\
Linf\_hier13x13x13d              & 39                                                        & 1.05                                   & 1                                   \\
Linf\_hier13x13x13e              & 39                                                        & 1.05                                   & 1                                  \\
Linf\_jjtabeltest3               & 56                                                        & 1                                   & 1                                   \\
Linf\_nine5d                     & 22                                                        & 1                                   & 1                                   \\
Linf\_nine12                     & 26                                                        & 1                                   & 1                                   \\
Linf\_ninenew                    & 24                                                        & 1                                   & 1                                   \\
Linf\_osorio                     & 30                                                        & 1                                   & 1                                   \\
Linf\_table1                     & 50                                                        & 1.6                                   & 1.6                                   \\
Linf\_table3                     & 52                                                        & 1                                   & 1                                   \\
Linf\_table6                     & 42                                                        & 1.05                                   & 1.05                                  \\
Linf\_table8                     & 34                                                        & 1                                   & 1                                   \\
Linf\_two5in6               & 40                                                        & 1                                   & 1     
\\
\hline
\end{tabular*}
\end{table}
\hspace{5mm}

One of the key advantages of our proposed method is the fast convergence in the Conjugate Gradient method to solve linear systems iteratively,  which allows it to be competitive with others, regarding running time. Additionally, our method was able to solve linear programming problems declared that were infeasible or Unknown for the other approaches; this is evident in the results for the \textit{Linf\_table3} problem, where our method is the only dense approach that could obtain an optimal solution.

Furthermore, our method shows consistent performance across different problem sizes. This is an important attribute of a robust optimization method, as it allows for reliable predictions and results in a wider range of applications.
\newline 

In order to evaluate the stability of the proposed method, we modify the density criterion. It is important to clarify that this criterion must be adjusted for each problem and it is just a proof of concept. In Table \ref{newresults1} we expose the features of the problems with the modified dense criteria. Furthermore, it is worth mentioning that by redefining the criterion, we were able to include additional problems that were initially categorized as lacking dense columns.

\begin{table}[h]
\caption{Tested problems with modified dense criteria}
\label{newresults1}
\begin{tabular*}{\hsize}{@{}@{\extracolsep{\fill}}lrrrr@{}}
\hline
  LP problem & {$m$}  & {$n$}  & {Nnz A}  & {ndense}   \\
\hline
Linf\_hier13                            & 3176                  & 3265                  & 14859                        & 24                         \\
Linf\_hier13x13x13a & 2992                  & 3096                  & 13293                        & 24   \\
Linf\_hier13x13x13b & 2992                  & 3096                  & 13293                     & 24                         \\
Linf\_hier13x13x13c                      & 2992                  & 3096                  & 13293                       & 24                         \\
Linf\_hier13x13x13d                      & 2992                  & 3096                  & 13293                      & 24                         \\
Linf\_hier13x13x13e                      & 2992                  & 3096                  & 13293                      & 22                         \\
Linf\_hier13x7x7d                        & 722                   & 825                  & 3404                        & 3                          \\
Linf\_jjtabeltest3                       & 2563                  & 2554                  & 9770                        & 28                         \\
Linf\_five20b                            & 48143     &50690                 & 239935                             & 14                         \\
Linf\_five20c                            & 47128                 & 48949                 & 233628                       & 36                         \\
Linf\_bts4                               & 55141                 & 63808                 & 265131                       & 20                         \\
Linf\_cbs                                & 3184                  & 3354                  & 12652                        & 2                          \\
Linf\_nine5d                             & 9853                  & 10532                 & 47308                        & 59                         \\
Linf\_nine12                             & 15815                 & 17592                 & 80293                        & 250                        \\
Linf\_ninenew                            & 9139                  & 10072                 & 46443                        & 27                         \\
Linf\_table1                             & 1308                  & 2734                  & 814                         & 2                          \\
Linf\_table3                             & 3477                  & 4165                  & 18846                        & 2                          \\
Linf\_table6                             & 1308                  & 1775                  & 7197                         & 2                          \\
Linf\_table8                             & 1342                  & 2539                  & 8800                        & 1                          \\
Linf\_two5in6                            & 5612                  & 5989                  & 27129                        & 21                         \\
Linf\_osorio                             & 10199                 & 19995                 & 69779                        & 2                          \\
fit2p                                    & 3000                  & 13525                 & 50284                       & 25                         \\
fit1p                                    & 627                   & 1677                  & 9868                        & 24                         \\
Linf\_dale                               & 15335                 & 24941                 & 89223                       & 2                          \\
Linf\_hier16                             & 5816                  & 6095                  & 28713                        & 9                          \\
Linf\_hier16x16x16a                      & 5564                  & 5851                  & 27102                        & 11                         \\
Linf\_hier16x16x16b                      & 5564                  & 5851                  & 27102                        & 11                         \\
Linf\_hier16x16x16c                      & 5564                  & 5851                  & 27102                       & 11                         \\
Linf\_hier16x16x16d                      & 5564                  & 5851                  & 27102                       & 11                         \\ 
\hline
\end{tabular*}
\end{table}

\hspace{4mm}

In Table \ref{newresults2} we summarize the results of problems related to running time, total iterations of the Interior-Point method and in the final column, ``Mean CG iter/IPM iter'' represents the mean number of Conjugate Gradient method iterations per Interior-Point method iteration, rounded to the nearest integer.

\begin{table}[h]
\caption{Results for modified dense criteria}
\label{newresults2}
\begin{tabular*}{\hsize}{@{}@{\extracolsep{\fill}}lrrr@{}}
\hline
  LP problem & {IPM Iterations}  & {Time(s)}  & {Mean CG Iter/ IPM Iter}   \\
\hline
Linf\_hier13        & 15                       & 0.46                     & 25                                         \\
Linf\_hier13x13x13a & 26                       & 0.02                     & 30                                         \\
Linf\_hier13x13x13b & 26                       & 0.02                     & 30                                         \\
Linf\_hier13x13x13c & 26                       & 0.68                     & 25                                         \\
Linf\_hier13x13x13d & 23                       & 0.61                     & 27                                         \\
Linf\_hier13x13x13e & 25                       & 0.64                     & 25                                         \\
Linf\_hier13x7x7d   & 20                       & 0.02                     & 3                                          \\
Linf\_jjtabeltest3  & 31                       & 0.58                     & 22                                         \\
Linf\_five20b       & 21                       & 158.01                   & 16                                         \\
Linf\_five20c       & 23                       & 867.45                   & 27                                         \\
Linf\_bts4          & 45                       & 19.64                    & 22                                         \\
Linf\_cbs           & 21                       & 0.05                     & 2                                          \\
Linf\_nine5d        & 17                       & 3.93                     & 15                                         \\
Linf\_nine12        & 15                       & 60.35                    & 65                                         \\
Linf\_ninenew       & 15                       & 8.87                     & 17                                         \\
Linf\_table1        & 23                       & 0.03                     & 2                                          \\
Linf\_table3        & 26                       & 0.27                     & 1                                          \\
Linf\_table6        & 25                       & 0.03                     & 2                                          \\
Linf\_table8        & 17                       & 0.02                     & 1                                          \\
Linf\_two5in6       & 20                       & 1.0                      & 22                                         \\
Linf\_osorio        & 13                       & 0.13                     & 1                                          \\
fit2p               & 19                       & 0.4                      & 25                                         \\
fit1p               & 19                       & 0.06                     & 22                                         \\
Linf\_dale          & 21                       & 0.29                     & 1                                          \\
Linf\_hier16        & 14                       & 4.53                     & 8                                          \\
Linf\_hier16x16x16a & 19                       & 4.73                     & 12                                         \\
Linf\_hier16x16x16b & 19                       & 4.74                     & 12                                         \\
Linf\_hier16x16x16c & 19                       & 4.74                     & 12                                         \\
Linf\_hier16x16x16d & 19                       & 4.88                     & 11   \\
\hline
\end{tabular*}
\end{table}

 For the results, we observe that, for all the cases, the mean number of iterations in the Conjugate Gradient method is positively correlated with the number of dense columns. This correlation arises because the matrix $W$ being solved involves the normal matrix of the dense part (Step $6$ of Algorithm \ref{newprecondperm}). Notably, in problems like \textit{Linf\_nine5d}, a significant number of dense columns are obtained, allowing for a successful problem solution. 
 \newline 
 
In all cases, $F=0$ was achieved, obtaining theoretical properties such as the minimum eigenvalue of $W$ (See Section \ref{propertiesofW}). This may occur because the number of sparse columns continues to be greater than or equal to the number of rows, since $D$ contains few columns compared to the original size of the constraint matrix $A$. A comparison with the initial criterion reveals an increase in the number of iterations and time, which is reasonable considering there are more dense columns. 
\newline 

 Furthermore, it was noted that, for the problems \textit{Linf\_five20b}, \textit{Linf\_five20c}, \textit{Linf\_bts4}, \textit{Linf\_cbs}, \textit{Linf\_ninenew}, \textit{Linf\_osorio} ($28\%$ of the initial $22$ tested problems), the running time required to find an optimal solution was reduced despite an increase in dense columns. This phenomenon can be attributed to the fact that the Cholesky factorization, which needs to be computed for each iteration of the Interior-Point method, originates from a matrix $S$ with a smaller number of columns and greater sparsity. Consequently, this results in a sparser Cholesky factorization process. This result also makes room for more sophisticated ways to define dense columns in the future.
 \newline 
 
\begin{figure}[H]
\centering
\includegraphics[width=.65\textwidth]{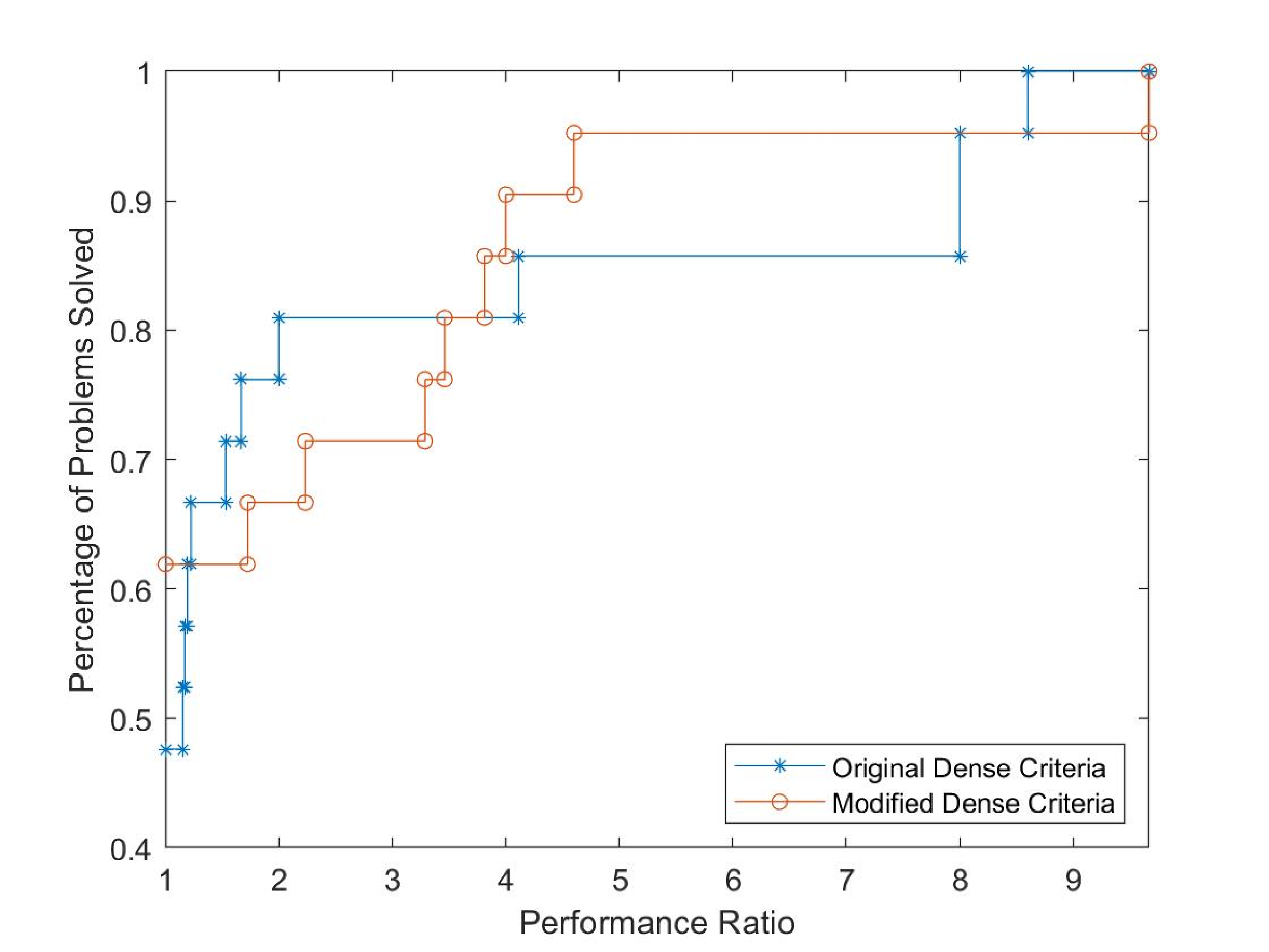}
\caption{ {\small Performance profile for the running time on a subset of problems. }}
\label{figura01}
\end{figure}

 Figure \ref{figura01} illustrates the performance profile regarding the running time required to solve the 21 common problems in Tables \ref{runningtime} and \ref{newresults2}. These problems are solved with the proposed preconditioner Chapter \ref{sec3} with PCx\_mod. The analysis focuses on examining the impact of varying the number of dense columns on the solving time. 
\newline 

Approximately $60\%$ of the tested problem, the running time is notably faster when dealing with a higher number of dense columns. This observation exposes the robustness of our approach, indicating that it can effectively adapt and perform well under different problem conditions, exposing its versatility and reliability in solving linear programming problems with a large number of dense columns.

\section{Conclusions}
\label{conclusion}
The omission of dense columns significantly extends the running time compared to approaches that account for them, as evident in our computational experiments. This paper introduces a novel approach for handling linear programming problems with dense columns in the constraint matrix. Theoretical demonstrations show that the final system to solve remains uniformly bounded as it approaches an optimal solution using Interior-Point methods.

A thorough comparison of our proposed preconditioner is conducted against existing approaches. In the realm of computational experiments, similar outcomes are observed in terms of the number of iterations and running time for $45\%$ of cases. Our approach consistently solved a higher number of linear programming problems, showing its robustness regarding running time and achieving the best performance in $37\%$ of tested problems.

In summary, the results from our experiments underscore the effectiveness and robustness of our proposed method for solving linear programming problems with dense columns. We trust that our method could prove invaluable for applications in fields such as engineering, finance, and operations research, where linear programming problems with dense columns are frequently encountered.

\section*{Acknowledgments}
This study was partially supported by the Conselho Nacional de Desenvolvimento Científico e Tecnológico -Brazil (CNPq) - Grants 141611/2020-0 and 313258/2021-0, and by the Coordenação de Aperfeiçoamento de Pessoal de Nível Superior - Brasil (CAPES) - Finance Code 001.


 \bibliographystyle{abbrv}
\bibliography{ref}

\end{document}